\documentclass[11pt]{amsart}

\usepackage[latin1]{inputenc}
\usepackage{amsmath,amsfonts,amssymb}
\usepackage{amscd}
\usepackage{graphicx}

\newcommand{\Hom}{\operatorname{Hom}}

\newcommand{\Ann}{\operatorname{Ann}}

\newtheorem{theorem}{Theorem}[section]
\newtheorem{proposition}[theorem]{Proposition}
\newtheorem{corollary}[theorem]{Corollary}
\newtheorem{lemma}[theorem]{Lemma}

\newtheorem{remark}[theorem]{Remark}

\begin{document}

\title[Generalized Reynolds ideals for non-symmetric algebras]
{Generalized Reynolds ideals for non-symmetric algebras}

\author[C. Bessenrodt, T. Holm, A. Zimmermann]
{Christine Bessenrodt \and Thorsten Holm \and Alexander Zimmermann}

\address{~~\newline
Christine Bessenrodt\newline
Institut f\"ur Algebra, Zahlentheorie und Diskrete Mathematik,
Universit\"at Hannover,
Welfengarten 1,
30167 Hannover, Germany}
\email{bessen@math.uni-hannover.de}

\address{~~\newline
Thorsten Holm\newline
Institut f\"ur Algebra und Geometrie,
Otto-von-Guericke-Universit\"at Magdeburg,
Postfach 4120, 39016 Magdeburg, Germany
\newline
and\newline
Department of Pure Mathematics,
University of Leeds,
Leeds LS2 9JT,
UK
}
\email{thorsten.holm@mathematik.uni-magdeburg.de}

\address{~~\newline
Alexander Zimmermann\newline Facult\'e de Math\'ematiques et LAMFA
(UMR 6140 du CNRS) , Universit\'{e} de Picardie, 33 rue St Leu,
80039 Amiens CEDEX 1, France}
\email{alexander.zimmermann@u-picardie.fr}


\bigskip

\begin{abstract}
We show how to extend the theory of generalized Reynolds ideals,
as introduced by B.\,K\"ulshammer, from symmetric algebras to
arbitrary finite-dimensional algebras (in positive characteristic).
This provides new invariants of the derived module categories
of finite-dimensional algebras, extending recent results
in \cite{Z} for symmetric algebras.
\end{abstract}

\maketitle

\section{Introduction}

Let $\Lambda$ be a finite-dimensional algebra over a field $k$
of characteristic $p>0$. The commutator space $K(\Lambda)$ is the
$k$-vector space
generated by all $[a,b]:=ab-ba$ where $a,b\in\Lambda$.
For any $n\ge 0$ set
$T_n(\Lambda) := \{x\in \Lambda\,\mid\,x^{p^n}\in K(\Lambda)\}$.
In \cite{Ku1}, B.\, K\"ulshammer defined  for any symmetric $k$-algebra
$\Lambda$ the generalized Reynolds ideals $T_n\Lambda^{\perp}$
as the orthogonal spaces with respect to the symmetrizing form
on the symmetric algebra $\Lambda$.
For an arbitrary finite-dimensional algebra $\Lambda$, its trivial extension
${\bf T}(\Lambda)$ is
a symmetric algebra.
For further details, definitions and background on generalized
Reynolds ideals and on trivial extensions
we refer to Sections~\ref{sec-Reynolds}
and~\ref{sec-trivialext} below.
\smallskip

The following are the main results of this note. They extend the
main result of \cite{Z} from symmetric algebras to arbitrary
finite-dimensional algebras.

\begin{theorem} Let $\Lambda$ be a finite-dimensional algebra
over a field of characteristic $p>0$. Let
${\bf T}(\Lambda)=\Lambda\ltimes \Lambda^*$ be the trivial extension.
Then the series of generalized Reynolds ideals
$T_i({\bf T}(\Lambda))^{\perp}$ of the trivial extension takes the form
$$Z(\Lambda)\ltimes \Ann_{\Lambda^*}(K( \Lambda))
\supseteq 0 \ltimes \Ann_{\Lambda^*}(T_1\Lambda)
\supseteq 0 \ltimes \Ann_{\Lambda^*}(T_2\Lambda)
\supseteq \ldots
$$
\end{theorem}

By a result of Rickard \cite[Theorem~3.1]{Rickard},
derived equivalent algebras have derived equivalent
trivial extension algebras.
As a consequence of this fact and our theorem above
we obtain some new numerical invariants of the
derived module categories of finite-dimensional algebras in
positive characteristic.

\begin{corollary} \label{cor-trivext}
Let $\Lambda$ be a finite-dimensional algebra over a perfect field
of characteristic $p>0$.
Then, for any $n\ge 0$, the codimension $\dim\Lambda - \dim T_n(\Lambda)$
is invariant under derived equivalences.
\end{corollary}

Note that the case $i=0$ gives the dimension
of $\Lambda/K( \Lambda)$ which is isomorphic to the degree zero
Hochschild homology; the latter is well-known to be invariant
under derived equivalences.
But for $i>0$ the above codimensions seem to be new derived
invariants, and do not seem to have an interpretation in terms
of other well-known invariants.
\smallskip

We can even prove the following more general and precise statement,
involving the structure as modules over the center.

\begin{corollary}\label{centremodules}
Let $\Lambda$ and $\Gamma$ be finite dimensional $k$-algebras over
a perfect field
$k$ of characteristic $p>0$. If the bounded derived categories of
$\Lambda$ and of $\Gamma$ are equivalent as triangulated categories
then the isomorphism induced between $Z(\Lambda)$ and $Z(\Gamma)$
by a twosided tilting complex
induces an isomorphism of the sequence of $Z(\Lambda)$-modules
$$\Ann_{\Lambda^*}(T_1(\Lambda))\supseteq \Ann_{\Lambda^*}(T_2(\Lambda))
\supseteq \Ann_{\Lambda^*}(T_3(\Lambda))\supseteq\cdots$$
and the sequence of $Z(\Gamma)$-modules
$$\Ann_{\Gamma^*}(T_1(\Gamma))\supseteq \Ann_{\Gamma^*}(T_2(\Gamma))
\supseteq \Ann_{\Gamma^*}(T_3(\Gamma))\supseteq\cdots$$
\end{corollary}

The paper is organized as follows. In Section \ref{sec-Reynolds}
we recall the definition and some known facts on generalized
Reynolds ideals for symmetric algebras. In Section \ref{sec-trivialext}
we collect several properties of trivial extensions which are
crucial for our purposes. Section \ref{sec-proof} is the core
part of this note, containing the proofs of our main results.

\section{Generalized Reynolds ideals} \label{sec-Reynolds}

In this section we briefly recall the definition of the sequence
of generalized Reynolds ideals as introduced by B. K\"ulshammer
\cite{Ku1}, \cite{Ku2}. For interesting recent developments on
this invariant we also refer to \cite{BHHKM}, \cite{HHKM},
\cite{Z}, \cite{Z2}.

Let $k$ a perfect field of characteristic $p > 0$.
Let $A$ be a finite-dimensional symmetric $k$-algebra with
associative, symmetric, nondegenerate $k$-bilinear form
$\langle -,-\rangle : A
\times A \rightarrow k$. For any subspace $M$ of $A$ we denote by
$M^{\bot}$ the orthogonal space of $M$ in $A$ with
respect to the form $\langle -,-\rangle$.
Moreover, let $K(A)$ be the $k$-subspace
of $A$ generated by all commutators $[a,b]:=a b - b a$,
$a, b \in A$. For any $n \geq 0$ set
$$T_n (A) = \left\{ x \in A \mid x^{p^n} \in K(A)\right\}.$$
Then, by \cite{Ku1}, for any $n\ge 0$, the orthogonal
space $T_n (A)^{\bot}$ is an ideal of the center
$Z(A)$ of $A$. These ideals are called {\it generalized Reynolds ideals}.
They form a descending sequence
$$Z(A) =  K(A)^{\perp} = T_0(A)^{\perp} \supseteq T_1(A)^{\perp} \supseteq T_2(A)^{\perp}
\supseteq T_3(A)^{\perp} \supseteq \ldots$$
In fact, B.\,K\"ulshammer showed in \cite{Ku1,Ku2} that there is
a mapping $\xi_n :
Z(A) \rightarrow Z(A)$ so that the equation
$$\langle\xi_n(z),x\rangle^{p^n} = \langle z,x^{p^n}\rangle$$
holds for any $z \in Z(A)$ and $x\in A/K(A)$. Moreover, he proved
that $\xi_n(Z(A)) = T_n (A)^{\bot}$.
Dually, for every $n\geq 1$, B.\,K\"ulshammer proved in \cite{Ku1,Ku2} that the
equation
$$\langle z,\kappa_n(x)\rangle^{p^n} = \langle z^{p^n},x\rangle$$
for any $z \in Z(A)$ and $x\in A/K(A)$ defines a mapping $\kappa_n
: A/K(A) \rightarrow A/K(A)$.  Considering the kernel and the image of these
maps leads to the following invariants,
$$\ker\kappa_n=\{\,x\in A/K(A)\,\mid\;\langle
z^{p^n},x\rangle=0\; \forall z\in Z(A)\,\} =: P_n(Z(A))^\perp/K(A)$$ and
\begin{eqnarray*}
\kappa_n(A/K(A))&=& T_n(Z(A))^\perp/K(A)\\
&=&\{\,x\in A/K(A)\,\mid\;\forall\; z\in Z(A):
z^{p^n}=0\Rightarrow \langle z,x\rangle=0\,\}.
\end{eqnarray*}

In \cite{HHKM} it has been shown that the sequence of generalized
Reynolds ideals is invariant under Morita equivalences. More
generally, the following theorem has been proven recently by
the third author.

\begin{proposition}[\cite{Z}, Theorem 1; \cite{Z2} Proposition 2.3
and Corollary 2.4]
\label{prop:zimmermann} %
Let $A$ and $B$ be finite-dimensional symmetric algebras over a perfect
field of positive characteristic $p$.
If $A$ and $B$ are derived equivalent, then
\begin{itemize}
\item[{(i)}]
there is an isomorphism $\varphi : Z(A) \rightarrow Z(B)$ between
the centers of $A$ and $B$ such that $\varphi(T_n (A)^{\bot}) = T_n
(B)^{\bot}$ for all positive integers $n$;
\item[{(ii)}]
there is an isomorphism $\psi: A/K(A)\rightarrow B/K(B)$ so that
$$P_n(Z(A))^\perp/K(A)\mbox{ is mapped by }\psi\mbox{ to }
P_n(Z(B))^\perp/K(B)$$ and so that
$$T_n(Z(A))^\perp/K(A)\mbox{ is
mapped by }\psi\mbox{ to }T_n(Z(B))^\perp/K(B).$$
\end{itemize}
\end{proposition}

We note that in the proof of \cite[Theorem 1]{Z} the
fact that $k$ is algebraically closed is never used.
The assumption on the field $k$ to be perfect is sufficient.
Hence the sequence of generalized Reynolds ideals gives a new
derived invariant for symmetric algebras over perfect
fields of positive characteristic.

The aim of this note is to show how one can extend this result
from symmetric to arbitrary finite-dimensional algebras, using
trivial extensions.

\section{Trivial extensions} \label{sec-trivialext}

Let $\Lambda$ be a finite-dimensional algebra over a field $k$.
We denote by $\Lambda^*$ the $k$-linear dual $\Hom_k(\Lambda,k)$
which becomes a $\Lambda$-$\Lambda$-bimodule by setting
$(a\varphi)(b)=\varphi(ba)$ and $(\varphi a)(b)=\varphi(ab)$ for all
$a,b\in\Lambda$ and all $\varphi\in\Lambda^*$.

The {\em trivial extension} ${\bf T}(\Lambda):= \Lambda \ltimes \Lambda^*$
is the $k$-algebra defined by the multiplication
$$(a,\varphi)\cdot (b,\psi) := (ab, a\psi+\varphi b) ~~~~~~~
\mbox{~~~for all}~~a,b\in\Lambda,~\varphi,\psi\in \Lambda^*.$$
Recall that an algebra $A$
is symmetric if there exists a $k$-linear map
$\pi:A\to k$ such that $\pi(ab)=\pi(ba)$ for all $a,b\in A$, and
such that the kernel of $\pi$ does not contain any nonzero
left or right ideals of~$A$.
The corresponding associative non-degenerate symmetric
$k$-bilinear form $\langle -,-\rangle$ on~$A$ is then given by
$$\langle a,b\rangle =\pi(ab)\; \mbox{~~~~for } a,b\in A\:.$$

Then we have the following well-known fact, which is the crucial
property of trivial extensions in our context.

\begin{proposition}
The trivial extension ${\bf T}(\Lambda)$ is a symmetric algebra, with
respect to the map $\pi:{\bf T}(\Lambda)\to k$, $(a,\varphi)\mapsto
\varphi(1)$.
\end{proposition}

\proof Clearly, for any $a,b\in\Lambda$ and $\varphi,\psi\in
\Lambda^*$ we have
$$\pi( (a,\varphi)\cdot (b,\psi)) = \psi(a)+\varphi(b)
=\pi( (b,\psi)\cdot (a,\varphi) ) .$$
Now let $I\lhd {\bf T}(\Lambda)$ be a left ideal contained in
the kernel of $\pi$. Let
$(b,\psi)\in I$. Since $I$ is a left ideal of ${\bf T}(\Lambda)$,
we get for all $a\in \Lambda$ and all $\varphi\in \Lambda^*$
that
$$(a,\varphi)\cdot (b,\psi)=(ab, a\psi+\varphi b)\in \ker\pi,$$
i.e., $\psi(a)+\varphi(b)=0$ for all $a\in\Lambda$, $\varphi\in
\Lambda^*$. In particular, by setting $\varphi=0$ we conclude that
$\psi(a)=0$ for all $a\in\Lambda$, and hence $\psi=0$. Then, also
$\varphi(b)=0$ for all
$\varphi\in\Lambda^*$, and hence $b=0$ as well. Therefore, $I=0$.

Similarly, $\ker\pi$ does not contain a nonzero right ideal.
\qed

\bigskip

We now collect some easy fundamental properties of trivial
extensions which will be used later in the proof of the main result.

For any subspace $V\le \Lambda$ we set
$$\Ann_{\Lambda^*}(V):= \{\varphi\in\Lambda^*\,\mid\,
\varphi(V)=0\},$$
i.e., those linear maps which vanish on $V$.

\begin{proposition} \label{prop-center}
Let $\Lambda$ be a finite-dimensional algebra
over a field $k$, with center $Z(\Lambda)$ and commutator subspace
$K( \Lambda)$.
Then the center of the trivial extension ${\bf T}(\Lambda)=\Lambda
\ltimes \Lambda^*$ has the form
$$Z({\bf T}(\Lambda)) = Z(\Lambda) \ltimes \Ann_{\Lambda^*}(K( \Lambda)).$$
\end{proposition}

\proof Let $(a,\varphi)\in Z({\bf T}(\Lambda))$. Clearly, $a\in Z(\Lambda)$.
Moreover, we then have
$$a\psi+\varphi b = b\varphi + \psi a \mbox{~~~~~for all~~}b\in\Lambda,~
\psi\in \Lambda^*.$$
The latter equality holds if and only if for all $c\in\Lambda$ we have
$$\psi(ca)+ \varphi(bc)=\varphi(cb) + \psi(ac).$$
Since $a\in Z(\Lambda)$, this holds precisely when
$\varphi(bc) = \varphi(cb)$ for all $b,c\in\Lambda$, i.e.,
when $\varphi(K( \Lambda))=0$, as claimed.
\qed

\begin{proposition} \label{prop-Ti}
Let $\Lambda$ be a finite-dimensional algebra
over a field $k$ of characteristic $p>0$, and let
${\bf T}(\Lambda)=\Lambda \ltimes \Lambda^*$ be the trivial extension.
Then for any $n\ge 0$ we get
$$T_n({\bf T}(\Lambda)) = T_n(\Lambda) \ltimes \Lambda^*.$$
\end{proposition}

\proof Recall that $T_n({\bf T}(\Lambda))= \{(a,\varphi)\in
{\bf T}(\Lambda)\,\mid\,(a,\varphi)^{p^n}\in K({\bf T}(\Lambda))\}$.
We have $(a,\varphi)^{p^n} = ((a,0)+(0,\varphi))^{p^n}$
which modulo the commutator subspace $K({\bf T}(\Lambda))$ becomes
$$(a,\varphi)^{p^n} \equiv (a,0)^{p^n} + (0,\varphi)^{p^n}
= (a^{p^n},0)$$
using that $(0\ltimes\Lambda^*)^2=0$.
In particular, $(a,\varphi)^{p^n}$ is in $K({\bf T}(\Lambda))$
if and only if $a^{p^n}\in K( \Lambda)$, i.e. $a\in T_n(\Lambda)$,
as claimed.
\qed
\bigskip

We need a final preparation for the main result,
namely an explicit description for the commutator
subspace $K({\bf T}(\Lambda))$ of a trivial extension.

\begin{lemma} \label{lemma-commutator}
Let $\Lambda$ be a finite-dimensional algebra
over a field. Then the following holds.
\begin{enumerate}
\item[{(1)}]
$K({\bf T}(\Lambda))=K(\Lambda) \ltimes [\Lambda,\Lambda^*]$
where $[\Lambda,\Lambda^*]$ denotes the commutator subspace of
$\Lambda^*$ of the algebra
$\Lambda$ acting on the bimodule $\Lambda^*$.
\item[{(2)}] Assume that $\Lambda$ is symmetric,
with symmetrizing form $\langle-,-\rangle$.
\begin{enumerate}
\item[{(a)}] The space $[\Lambda,\Lambda^*]$ is generated by all
linear maps $$\varphi_{[a,b]}:=\langle -,ab-ba\rangle~~~~
\mbox{where $a,b\in\Lambda$}.$$
\item[{(b)}] We have $[\Lambda,\Lambda^*] = \Ann_{\Lambda^*}(Z(\Lambda)).$
\end{enumerate}
\end{enumerate}
\end{lemma}

\proof
(1) First, for any $a,b\in \Lambda$ we have
$$(a,0)\cdot (b,0)-(b,0)\cdot (a,0)=(ab-ba,0),$$
i.e., this generates $K(\Lambda)\ltimes 0$. Secondly,
for any $a\in\Lambda$ and $\varphi\in\Lambda^*$ we get
$$(a,0)\cdot (0,\varphi)-(0,\varphi)\cdot(a,0)=(0,a\varphi-\varphi a)$$
which generates $0\ltimes [\Lambda,\Lambda^*]$.
This shows the inclusion '$\supseteq$'.

Conversely, an arbitrary generator of $K({\bf T}(\Lambda))$ has the form
$$[(a,\varphi),(b,\psi)] = (ab,a\psi+\varphi b) - (ba,b\varphi+\psi a)
= (ab-ba,a\psi-\psi a + \varphi b-b\varphi)
$$
which clearly is in $K(\Lambda)\ltimes [\Lambda,\Lambda^*]$.
\smallskip

(2\,a) The commutator space $[\Lambda,\Lambda^*]$ is generated
by all $a\varphi - \varphi a$ where $a\in\Lambda$, and $\varphi\in
\Lambda^*$. Since $\Lambda$ is finite-dimensional, the elements
of $\Lambda^*$ are of the form $\varphi=\varphi_b := \langle -,b\rangle$
where $b\in \Lambda$. It suffices to show that
$a\varphi_b - \varphi_b a = \varphi_{ab-ba}$. This holds, since for every
$x\in\Lambda$ we have
\begin{eqnarray*}
(a\varphi_b - \varphi_b a)(x)&=& \varphi_b(xa) - \varphi_b(ax) =
\langle xa-ax,b\rangle = \langle x,ab\rangle - \langle b,ax\rangle \\
&=& \langle x,ab\rangle - \langle ba,x\rangle =
\langle x,ab\rangle - \langle x,ba\rangle \\
&=& \varphi_{ab-ba}(x).\\
\end{eqnarray*}

\vskip-0.45cm

(2\,b) Let $a\in\Lambda$, $\varphi\in
\Lambda^*$ and $z\in Z(\Lambda)$. Then we have
$$(a\varphi -\varphi a)(z)=\varphi(az)-\varphi(za)=\varphi(az-za)
=\varphi(0)=0,$$
i.e., $[\Lambda,\Lambda^*]\subseteq \Ann_{\Lambda^*}(Z(\Lambda))$.

Conversely, let $\varphi\in \Ann_{\Lambda^*}(Z(\Lambda))$.
Since $\Lambda$ is finite-dimensional,
there exists an element $b\in\Lambda$ such that $\varphi = \varphi_b
=\langle -,b\rangle$. By assumption on $\varphi$ we have that
$$b\in Z(\Lambda)^{\perp} = (K(\Lambda)^{\perp})^{\perp} = K(\Lambda).$$
From part (2\,a) it now follows that $\varphi=\varphi_b\in [\Lambda,\Lambda^*]$.
\qed

\section{Proof of the main results} \label{sec-proof}

This section is devoted to proving the main results of this article.

The following main step shows that
indeed the generalized Reynolds ideals for the trivial extension
${\bf T}(\Lambda)$ are
closely related to the algebra $\Lambda$ itself.

\begin{theorem} \label{thm-reynolds}
Let $\Lambda$ be a finite-dimensional algebra
over a field of characteristic $p>0$, and let
${\bf T}(\Lambda)=\Lambda\ltimes\Lambda^*$ be its trivial extension.
\begin{enumerate}
\item[{(1)}] We have $T_0({\bf T}(\Lambda))^\perp=Z({\bf T}(\Lambda))=
Z(\Lambda)\ltimes \Ann_{\Lambda^*}(K(\Lambda))$.
\item[{(2)}] For all $n\ge 1$ the generalized Reynolds ideals of
the trivial extension are of the form
$T_n({\bf T}(\Lambda))^{\perp} = 0 \ltimes \Ann_{\Lambda^*}(T_n\Lambda).$
\item[{(3)}]
For all $n\ge 1$ we have\\
$T_n(Z({\bf T}(\Lambda)))^\perp/K({\bf T}(\Lambda))=
0\ltimes \left(\Ann_{\Lambda^*}(T_n(Z(\Lambda)))/[\Lambda,\Lambda^*]\right).$
\item[{(4)}] For all $n\ge 1$ we have\\
$P_n(Z({\bf T}(\Lambda)))^\perp/K({\bf T}(\Lambda))=
\Lambda/K(\Lambda)\ltimes
\left(\Ann_{\Lambda^*}(P_n(Z(\Lambda)))/[\Lambda,\Lambda^*]\right)$.
\end{enumerate}
\end{theorem}

\proof (1) The first equality follows directly from the
general fact that for a symmetric algebra the orthogonal
space of the commutator subspace is the center (see e.g. \cite{Ku2}).
The second equality is
Proposition \ref{prop-center}.
\smallskip

(2) Let us fix some $n\ge 1$. Recalling the definition
of the symmetric bilinear form
on ${\bf T}(\Lambda)$ and using Proposition \ref{prop-Ti} we have
$$T_n({\bf T}(\Lambda))^{\perp} = \{ (b,\psi)\in {\bf T}(\Lambda)\,\mid\,
\psi(a)+\varphi(b)=0\mbox{~~~for all~~}a\in T_n(\Lambda),~\varphi\in
\Lambda^*\}.$$
Setting $a=0\in T_n(\Lambda)$
we conclude that $b=0$ and then we get
$$T_n({\bf T}(\Lambda))^{\perp} = \{ (0,\psi)\in {\bf T}(\Lambda)\,\mid\,
\psi(T_n\Lambda)=0\},$$
as claimed.
\smallskip

(3) First note that, since $Z({\bf T}(\Lambda))$ is commutative,
the corresponding
commutator space $K(Z({\bf T}(\Lambda)))$ becomes zero.
In particular,
$$(a,\varphi)\in T_n(Z({\bf T}(\Lambda)))\Longleftrightarrow
(0,0) = (a,\varphi)^{p^n}= (a,0)^{p^n}+(0,\varphi)^{p^n} = (a^{p^n},0).$$
Hence,
$T_n(Z({\bf T}(\Lambda)))=T_n(Z(\Lambda))\ltimes \Ann_{\Lambda^*}(K(\Lambda))\;.$
Therefore, we can conclude for the orthogonal space that
$$(a,\varphi)\in T_n(Z({\bf T}(\Lambda)))^\perp~~\Longleftrightarrow
~~\langle (a,\varphi),(z,\psi)\rangle = \varphi(z)+\psi(a)=0$$
\mbox{ for all }
$(z,\psi)\in T_n(Z(\Lambda))\ltimes \Ann_{\Lambda^*}(K(\Lambda))$.
In particular, setting $z=0$ one gets $a\in K(\Lambda)$. Moreover,
choosing $\psi=0$ shows that $\varphi\in\Ann_{\Lambda^*}(T_n(Z(\Lambda)))$.
Now invoking Lemma \ref{lemma-commutator}
gives the required statement.
\smallskip

(4) Using the definition of the
symmetrizing bilinear form on ${\bf T}(\Lambda)$
and the fact that $0\ltimes \Lambda^*$ is a nilpotent ideal of square zero
we have that
$$(a,\varphi)\in P_n(Z({\bf T}(\Lambda)))^\perp/K({\bf T}(\Lambda))
\Longleftrightarrow
\varphi(z^{p^n})=0\mbox{~~~for all~} z\in Z(\Lambda).$$
Now recall that by definition
$P_n(Z(\Lambda))$ is the space generated by all $z^{p^n}$ where
$z\in Z(\Lambda)$. Then
we can use Lemma \ref{lemma-commutator}
to deduce the desired statement.
\qed
\medskip

Corollary \ref{cor-trivext}, as stated in the introduction, can
now be deduced from the above theorem as follows.
\medskip

{\em Proof of Corollary \ref{cor-trivext}.} Suppose $\Lambda$
and $\Gamma$ are derived equivalent
finite-dimensional algebras over a perfect field
of positive characteristic. Then we consider their trivial
extensions ${\bf T}(\Lambda)$ and ${\bf T}(\Gamma)$, respectively. These
are derived equivalent, by a result of J. Rickard
\cite{Rickard}. Moreover, the trivial extensions are symmetric algebras.
Hence, they have generalized Reynolds ideals, and by Theorem
\ref{thm-reynolds} the sequences of these take the form
$$Z(\Lambda)\ltimes \Ann_{\Lambda^*}(K( \Lambda))
\supseteq 0 \ltimes \Ann_{\Lambda^*}(T_1\Lambda)
\supseteq 0 \ltimes \Ann_{\Lambda^*}(T_2\Lambda)
\supseteq \ldots
$$
and analogously for ${\bf T}(\Gamma)$
$$Z(\Gamma)\ltimes \Ann_{\Gamma^*}(K( \Gamma))
\supseteq 0 \ltimes \Ann_{\Gamma^*}(T_1\Gamma)
\supseteq 0 \ltimes \Ann_{\Gamma^*}(T_2\Gamma)
\supseteq \ldots
$$
According to \cite{Z}, these sequences are invariants of the
derived equivalent symmetric algebras ${\bf T}(\Lambda)$ and ${\bf T}(\Gamma)$.
In particular, the dimensions in each step have to coincide.
Recall that derived equivalent algebras have isomorphic centers, so
$Z(\Lambda)\cong Z(\Gamma)$. Now we compare the remaining
dimensions in the
above sequences. Note that for any $n\ge 0$ we have
$\dim\Ann_{\Lambda^*}(T_n\Lambda)=\dim\Lambda - \dim T_n\Lambda$
and similarly for $\Gamma$. Since these dimension have to agree,
we get that the codimensions $\dim\Lambda - \dim T_n\Lambda$
are derived invariant, as claimed.
\qed

\medskip

\begin{remark} \label{remnothingnew}
{\em
We are going to explain that our main result,
Theorem \ref{thm-reynolds} above,
really extends the corresponding result from \cite{Z}
for symmetric algebras. In fact, let us start with a finite-dimensional
symmetric algebra
$\Lambda$. Then the nondegenerate symmetrizing form
$\langle -,-\rangle$ on
$\Lambda$ induces an isomorphism $\lambda:\Lambda\longrightarrow
\Lambda^*$ as $\Lambda$-$\Lambda$-bimodules by setting
$$\lambda(b)=\left(\,a\mapsto \langle a,b\rangle\,\right).$$
Moreover, $K(\Lambda)$, and hence also
$T_n(\Lambda)$, is a $Z(\Lambda)$-submodule of $\Lambda$.
We now have two series of Reynolds ideals to compare. Namely, the one
for the symmetric algebra $\Lambda$ itself
$$Z(\Lambda)\supseteq T_1(\Lambda)^{\perp}\supseteq T_2(\Lambda)^{\perp}
\supseteq \ldots
$$
and the one for the trivial extension ${\bf T}(\Lambda)$ which,
by Theorem \ref{thm-reynolds}, takes the form
$$Z(\Lambda)\ltimes \Ann_{\Lambda^*}(K( \Lambda))
\supseteq 0 \ltimes \Ann_{\Lambda^*}(T_1(\Lambda))
\supseteq 0 \ltimes \Ann_{\Lambda^*}(T_2(\Lambda))
\supseteq \ldots
$$
We claim that the structure of the ideals occurring in these two series
are the same. In fact, for the former series, the $T_n(\Lambda)^{\perp}$
are ideals of $Z(\Lambda)$ by multiplication in the ring $Z(\Lambda)$.
For the latter series, take elements $(z,\psi)\in
Z(\Lambda)\ltimes \Ann_{\Lambda^*}(K( \Lambda))$ and
$(0,\varphi)\in 0 \ltimes \Ann_{\Lambda^*}(T_n(\Lambda))$ for some
$n\ge 1$. By the multiplication rule for the trivial extension we
have $(z,\psi)\cdot (0,\varphi) = (0,z\varphi)$ and similarly
$(0,\varphi)\cdot (z,\psi) = (0,\varphi z)$. Hence, the ideal
structure is given by the $Z(\Lambda)$-action on
$\Ann_{\Lambda^*}(T_n(\Lambda))$. Now our claim follows from the
observation that under the above identification,
$\Ann_{\Lambda^*}(T_n(\Lambda))\cong T_n(\Lambda)^{\perp}$,
$\langle -,t\rangle \leftrightarrow t$, the $Z(\Lambda)$-bimodule
action on $\Ann_{\Lambda^*}(T_n(\Lambda))$ corresponds precisely
to multiplication in $Z(\Lambda)$.
\smallskip

Now let us consider the invariants occurring in Proposition
\ref{prop:zimmermann}, part (ii), and Theorem \ref{thm-reynolds},
parts (3),(4), respectively. We first have to compare
$T_n(Z(\Lambda))^{\perp}/K(\Lambda)$ and
$0\ltimes \Ann_{\Lambda^*}(T_n(Z(\Lambda)))/[\Lambda,\Lambda^*]$.
Note that the former is contained in $\Lambda/K(\Lambda)$,
whereas the latter is contained in
${\bf T}(\Lambda)/K({\bf T}(\Lambda)).$
Again, we use the map $\lambda:\Lambda\to\Lambda^*$,
$b\mapsto \langle-,b\rangle$. Under this isomorphism,
$T_n(Z(\Lambda))^{\perp}$ corresponds to
$\Ann_{\Lambda^*}(T_n(Z(\Lambda)))$.
Moreover, $K(\Lambda)$ corresponds similarly to
$$\Ann_{\Lambda^*}(K(\Lambda)^{\perp}) = \Ann_{\Lambda^*}(Z(\Lambda))
= [\Lambda,\Lambda^*]$$
where the last equation is Lemma \ref{lemma-commutator}\,(2b).
\smallskip

The above argument applies analogously to the spaces
$P_n(Z(\Lambda))$ and then shows that, for $\Lambda$ symmetric,
$P_n(Z(\Lambda))^{\perp}/K(\Lambda)$
corresponds to the second argument in the invariant
$$P_n(Z({\bf T}(\Lambda)))^{\perp}/K({\bf T}(\Lambda))
 = \Lambda/K(\Lambda) \ltimes (\Ann_{\Lambda^*}(P_n(Z(\Lambda)))/
[\Lambda,\Lambda^*]).
$$
Taken together, the above considerations of this remark show that
for a symmetric algebra $\Lambda$ the invariants obtained from the
trivial extension ${\bf T}(\Lambda)$ carry the same information
as the invariants obtained from $\Lambda$ directly. In this sense,
our construction genuinely extends the results from \cite{Z}, \cite{Z2}
from symmetric to arbitrary finite-dimensional algebras.}
\end{remark}

%
%
%
%

\smallskip

As the final step we can now give the proof of our second corollary
from the introduction.
\medskip

{\em Proof of Corollary \ref{centremodules}.}
Let $\Lambda$ and $\Gamma$ be two derived equivalent $k$-algebras.
By Rickard's main theorem we know that the derived
equivalence can be replaced by a derived equivalence of standard
type, i.e. given by tensoring with a two-sided tilting complex.
Such an equivalence gives an isomorphism
$\zeta:Z({\bf T}(\Lambda))\rightarrow Z({\bf T}(\Gamma))$
and by Proposition~\ref{prop:zimmermann}
we know that then the decreasing sequence of $Z({\bf T}(\Lambda))$-ideals
$T_n({\bf T}(\Lambda))^\perp$ is mapped to the decreasing sequence
$T_n({\bf T}(\Gamma))^\perp$. But the Remark~\ref{remnothingnew}
implies that this is just the same as saying that
the isomorphism $\zeta$ induces a $Z(\Lambda)$-module isomorphism
$\Ann_{\Lambda^*}(K(\Lambda))\rightarrow \Ann_{\Gamma^*}(K(\Gamma))$
so that the sequence of submodules
$$\Ann_{\Lambda^*}(T_1\Lambda)\supseteq
\Ann_{\Lambda^*}(T_2\Lambda)\supseteq
\Ann_{\Lambda^*}(T_3\Lambda)\supseteq\dots$$ is
mapped to the sequence of submodules
$$\Ann_{\Gamma^*}(T_1\Gamma)\supseteq
\Ann_{\Gamma^*}(T_2\Gamma)\supseteq \Ann_{\Gamma^*}(T_3\Gamma)
\supseteq\dots\;.$$
This completes the proof of the corollary.\qed


\end{document}